\magnification=1200
\def\A{{\cal A}}
\def\B{{\cal B}}
\def\C{{\bf C}}
\def\F{{\cal F}}
\def\I{{\cal I}}
\def\K{{\cal K}}
\def\N{{\bf N}}

\def\U{{\cal U}}
\def\Z{{\bf Z}}
\def\hal{{\vrule height 10pt width 4pt depth 0pt}}

\centerline{\bf Consistency of a counterexample to Naimark's problem}
\medskip

\centerline{Charles Akemann and Nik Weaver\footnote{*}{Partially supported
by NSF grant DMS-0070634. The second author wishes to thank Gert
Pedersen for motivating his work on this problem.}}
\bigskip
\bigskip

{\narrower{
\noindent \it
We construct a C*-algebra that has only one irreducible representation
up to unitary equivalence but is not isomorphic to the algebra of compact
operators on any Hilbert space. This answers an old question of Naimark.
Our construction uses a combinatorial statement called the diamond principle,
which is known to be consistent with but not provable from the standard
axioms of set theory (assuming those axioms are consistent). We prove
that the statement ``there exists a counterexample to Naimark's problem
which is generated by $\aleph_1$ elements'' is undecidable in standard
set theory.
\bigskip}}

Let $\K(H)$ denote the C*-algebra of compact operators on a complex,
not necessarily separable, Hilbert space $H$. In [14] Naimark observed
that every irreducible representation (irrep) of $\K(H)$ is unitarily
equivalent to the identity representation, so that each of these algebras
has only one irrep up to equivalence, and in [15] he asked whether this
property characterizes the algebras $\K(H)$. In other words: if $\A$ is
a C*-algebra with only one irrep up to unitary equivalence, is $\A$
isomorphic to some $\K(H)$? We call any algebra $\A$ that satisfies the
premise of this question but not its conclusion a {\it counterexample to
Naimark's problem}.
\medskip

The problem was quickly settled in the separable case. Building on results
in [15], Rosenberg [17] showed that there are no separable (indeed, no
separably acting) counterexamples to Naimark's problem. Around the same time,
Kaplansky [10] introduced the so-called ``type I'' (or ``GCR'') C*-algebras
and began developing their representation theory. (See [3, 16, 20] for
general background on type I C*-algebras.) This development was
carried further by Fell and Dixmier, who showed in particular that any
two irreps of a type I C*-algebra with the same kernel are unitarily
equivalent [6] and conversely, any separable C*-algebra that is not type I
has inequivalent irreps with the same kernel [5]. As it is easy to see that
no type I C*-algebra can be a counterexample to Naimark's problem, the
latter result (partially) recovers Rosenberg's theorem, doing so,
moreover, in the context of a general theory.
\medskip

Next, in a celebrated paper Glimm [8] gave several characterizations of
separable type I C*-algebras, showing in particular that every separable
C*-algebra which is not type I has {\it uncountably many} inequivalent
irreps. Since counterexamples to Naimark's problem cannot be type I, this
re-established the nonexistence of separable counterexamples in an
especially dramatic way.
\medskip

Some of Glimm's results were initially proven without assuming separability,
and subsequent work by Sakai [18, 19] went further in this direction.
However, separability was never removed from the reverse implication
of the equivalence ``type I $\Leftrightarrow$ irreps with the same kernel
are equivalent,'' and it was recognized that Naimark's problem
potentially represented a fundamental obstruction to a nonseparable
generalization of this result. However, it was reasonable to expect that
there were no counterexamples to Naimark's problem because it seemed
likely that the separability assumption could be dropped in Glimm's
theorem that non type I C*-algebras have uncountably many inequivalent
irreps. But this was never achieved.
\medskip

That is the background for the present investigation. We construct a
(necessarily nonseparable and not type I) counterexample to Naimark's problem.
Our construction is not carried out in ZFC (Zermelo-Frankel set theory with
the axiom of choice): it requires Jensen's ``diamond'' principle, which
follows from G\"odel's axiom of constructibility and hence is relatively
consistent with ZFC. (See [9] or [13] for general background on set theory).
The diamond principle has been used to prove a variety of consistency results
in mainstream mathematics; see, e.g., [4, 21, 22].
\medskip

Presumably the existence of a counterexample to Naimark's problem is
independent of ZFC, but we have not yet been able to show this. We can
prove the relative consistency of the assertion ``no C*-algebra generated
by $\aleph_1$ elements is a counterexample to Naimark's problem''; indeed,
this follows easily from [8]. Since our counterexample is generated by
$\aleph_1$ elements, if follows that the existence of an
$\aleph_1$-generated counterexample is independent of ZFC (see
Corollary 7).
\medskip

The basic idea of our construction is to create a nested transfinite
sequence of separable C*-algebras $\A_\alpha$, for $\alpha < \aleph_1$,
each equipped with a distinguished pure state $f_\alpha$, such that for any
$\alpha < \beta$, $f_\beta$ is the unique extension of $f_\alpha$
to a state on $\A_\beta$. At the same time, for each $\alpha$ (or at
least ``enough'' $\alpha$) we want to
select a pure state $g_\alpha$ on $\A_\alpha$ which is not equivalent to
$f_\alpha$, such that $g_\alpha$ has a unique state extension $g_\alpha'$
to $\A_{\alpha + 1}$, and $g_\alpha'$ is equivalent to $f_{\alpha + 1}$.
Thus, we build up a continually expanding pool of equivalent pure states
in an attempt to ensure that all pure states on $\A = \bigcup \A_\alpha$ will
be equivalent.
\medskip

There are two significant difficulties with this approach, which incidentally
is essentially the only way an $\aleph_1$-generated counterexample to
Naimark's problem could be constructed. First, there is
the technical problem of finding a suitable algebra $\A_{\alpha + 1}$ which
contains $\A_\alpha$ and such that $f_\alpha$ and $g_\alpha$ extend uniquely
to $\A_{\alpha + 1}$ and become equivalent. We accomplish this via a crossed
product construction which takes advantage of a powerful recent result of
Kishimoto, Ozawa, and Sakai that, together with earlier work of Futamura,
Kataoka, and Kishimoto, ensures the existence of automorphisms on
separable C*-algebras which relate inequivalent pure states in a certain
manner. The second
fundamental challenge is to choose the states $g_\alpha$ in such a way that
{\it all} pure states are eventually made equivalent to one another. This
is especially troublesome because pure states can be expected to proliferate
exponentially as $\alpha$ increases (every pure state on $\A_\alpha$ extends
to at least one, but probably many, pure states on $\A_{\alpha + 1}$) and
we can only make a single pair of states equivalent at each step. We
handle this issue by using the diamond principle, one version of which states
that it is possible to select a single vertex from each level of the standard
tree of height and width $\aleph_1$, such that every path down the tree
contains in some sense ``many'' selected vertices. In our application
the vertices of the tree at level $\alpha$ model the states on $\A_\alpha$
and diamond informs us how to choose $g_\alpha$. Then every pure state
on $\A$ induces a path down the model tree, and hence its restriction to
some $\A_\alpha$ equals $g_\alpha$, so that all pure states are indeed
taken care of at some point in our construction.
\bigskip
\bigskip

\noindent {\bf 1. Unique extension of pure states to crossed products}
\bigskip

In this section we consider the problem: if $f$ and $g$ are inequivalent
pure states on a C*-algebra $\A$, is it possible to find a C*-algebra
$\B$ which contains $\A$ such that (1) $f$ and $g$ have unique state
extensions $f'$ and $g'$ to $\B$ and (2) $f'$ and $g'$ are equivalent?
We find that the answer is yes if $\A$ is simple, separable, and unital,
and moreover we can ensure that $\B$ is also simple, separable, and unital,
with the same unit as $\A$.
\medskip

A net $(a_\lambda)$ of positive, norm-one elements of a C*-algebra $\A$
is said to {\it excise} a state $f$ on $\A$ if $\|a_\lambda xa_\lambda -
f(x)a_\lambda^2\|
\to 0$ for all $x \in \A$ [1, p.\ 1239].  By [1, Proposition 2.2], for
each pure state $f$ of a $C^*$-algebra $A$, there is a decreasing net
$(a_\lambda)$ that excises $f$.  As described in [1, p.\ 1240; see also
16, p.\ 87], the support projection $p$ of $f$ is a rank 1 projection in
$\A^{**}$.  Further, as shown in the last 7 lines of the proof of [1,
Proposition 2.3], $a_\lambda \to p$ in the strong operator topology on
$A^{**}$.
\bigskip

\noindent {\bf Lemma 1.} {\it Suppose $f$ and $g$ are inequivalent
pure states on a C*-algebra $\A$, $(a_\lambda)$ excises $f$, and
$(b_\lambda)$ excises $g$. Also assume that $(b_\lambda)$ is decreasing.
Then $\|a_\lambda x b_\lambda\| \to 0$ for all $x \in \A$.}
\medskip

\noindent {\it Proof.} Note that we are assuming $(a_\lambda)$ and
$(b_\lambda)$ have the same index set. This is easy to arrange by
replacing two different index sets with their product.
Let $q$ denote the rank one projection in
$\A^{**}$ that supports the pure state $g$. Suppose the lemma
is false for some $x \in \A$.  Then, multiplying $x$ by a large enough
positive number and passing to a subnet, we can assume that
$\|a_\lambda x b_\lambda\| > 1$ for all $\lambda$.
\medskip

We first show that $f(xb_\kappa x^*)\searrow 0$.
In $\A^{**}$, $xb_\kappa x^*\searrow xqx^*$.  Let $y$ denote the
central cover of $q$ in $\A^{**}$.  Since $g$ is pure, $y\A^{**}$ is
isomorphic to some $B(H)$. Thus $f$ must be 0 on $y\A^{**}$
lest it be equivalent to $g$. Therefore
$$f(xb_\kappa x^*)\searrow f(xqx^*)=f(xyqx^*)=f(y(xqx^*))=0.$$

Fix $\kappa_0$ so that $f(xb_{\kappa_0} x^*)< 1/2$, and note that
because $(a_\lambda)$ excises $f$, there exists $\lambda_0$ such that
$\lambda\geq \lambda_0$ implies $\|a_\lambda xb_{\kappa_{0}} x^*a_\lambda -
f(xb_{\kappa_0} x^*)a_\lambda^2\|< 1/2$.  Choose $\lambda$ so that
$\lambda \ge \lambda_0$ and $\lambda \ge \kappa_0$.
Using $b_\lambda \le b_{\kappa_{0}}$, we have
$$\eqalign{\|a_\lambda xb_\lambda\|^2
&\leq \|a_\lambda x b_\lambda^{1/2}\|^2
= \|a_\lambda xb_\lambda x^*a_\lambda\|
\le \|a_\lambda xb_{\kappa_0} x^*a_\lambda\|\cr
&\le \|a_\lambda xb_{\kappa_0}x^*a_\lambda - f(xb_{\kappa_0}x^*)a_\lambda^2\|
+f(xb_{\kappa_0}x^*)\|a_\lambda^2\|\cr
&< 1/2 + 1/2 = 1,\cr}$$
a contradiction.\hfill\hal
\bigskip

Let $\theta$ be an action of a discrete group $G$ on a unital C*-algebra
$\A$. The reduced crossed product $\A \times_r G$ can be defined
as the C*-algebra which acts on the right Hilbert module $l^2(G;\A)$ and
is generated by (1) for each $x \in \A$, the multiplication operator
$M_x$ defined by $M_x\phi(g) = \theta_{g^{-1}}(x)\phi(g)$ and
(2) for each $h \in G$ the translation operator $T_h$ defined by
$T_h\phi(g) = \phi(h^{-1}g)$. Observe that the map $x \mapsto M_x$
isomorphically embeds $\A$ in $\A \times_r G$, so we can regard $\A$
as a subalgebra of $\A \times_r G$ by identifying $x$ with $M_x$.
Denote the identity of $G$ by $e$.
\bigskip

\noindent {\bf Theorem 2.} {\it Let $\A$ be a unital C*-algebra,
let $f$ be a pure state on $\A$, and let $\theta$ be an action of a discrete
group $G$ on $\A$. Embed $\A$ in $\A \times_r G$ in the manner just
described. Then $f$ has a unique state extension to $\A \times_r G$
if and only if $f$ is inequivalent to $f\circ\theta_g$ for all
$g \neq e$.}
\medskip

\noindent {\it Proof.} ($\Rightarrow$) Suppose that for some $g \in G$
not the identity, $f$ is equivalent to $f\circ\theta_g$;
we must show that $f$ does not extend uniquely to the crossed product.
Let $u \in \A$ be a unitary such that $f = u^*(f\circ\theta_g)u$
[16, Proposition 3.13.4].  Let $\pi: \A \to \B(H_f)$ be the GNS
representation associated to $f$ and let $\xi \in H_f$ be the image of the unit
of $\A$, so that $f(y) = \langle \pi(y)\xi, \xi\rangle$ for all $y \in \A$.
Let $b \in \A$ be any positive, norm-one element such that $f(b) = 1$;
then $\pi(b)(\xi) = \xi$. Furthermore,
$\theta_g(u^*bu)$ is also a positive, norm-one element such
that $f(\theta_g(u^*bu)) = f(b) = 1$, so the same is true of
this element, and it follows that $\pi(b\theta_g(u^*bu))(\xi) = \xi$.
Hence $\|b\theta_g(u^*bu)\| = 1$.
\medskip

Now let $x = \theta_g(u^*)$. Let $b \in \A$
be any positive, norm-one element such that $f(b) = 1$ and let $y \in \A$.
We will show that $\|y - b(xT_g)b\| \geq 1$ (computing the norm in the
crossed product), which implies non-unique extension by [2, Theorem 3.2].
Observe first that
$$b(xT_g)b = bx\theta_g(b)T_g = b\theta_g(u^*bu)\theta_{g^{-1}}(u)T_g = cT_g$$
where $c \in \A$ has norm $1$.  (It is the product of $b\theta_g(u^*bu)$, which
we showed above has norm $1$, with the unitary $\theta_{g^{-1}}(u)$.)
Now to estimate the norm of $y - b(xT_g)b = y - cT_g$, consider the element
$\delta_e$ of $l^2(G;\A)$ which satisfies $\delta_e(e) = 1_\A$ (the unit of
$\A$) and $\delta_e(h) = 0$ for all $h \neq e$. Since $g \neq e$, we have
$(y - cT_g)(\delta_e) = \phi$ where $\phi(e) = y$, $\phi(g) =
-\theta_{g^{-1}}(c)$, and $\phi(h) = 0$ for all other $h$. The norm of
$\delta_e$ in $l^2(G;\A)$ is
$\|\sum_h \delta_e(h)^*\delta_e(h)\|^{1/2}_\A = \|1_\A\|^{1/2}_\A = 1$
and the norm of $\phi$ is
$$\left\|\sum_h \phi(h)^*\phi(h)\right\|^{1/2}_\A
= \|y^*y + \theta_{g^{-1}}(c)^*\theta_{g^{-1}}(c)\|^{1/2}_\A \geq 1$$
since $\|\theta_{g^{-1}}(c)\| = \|c\| = 1$.  Thus the norm of $y - cT_g$
is at least 1, as we needed to show.
\medskip

($\Leftarrow$) Suppose $f$ is not equivalent to $f\circ\theta_g$ for
any $g \neq e$. To verify that $f$ has a unique extension to the crossed
product, according to [2, Theorem 3.2] we must, for every element $z$
of the crossed product and every $\epsilon > 0$, find $x \in \A$ and a
positive norm-one element $b \in \A$ such that $f(b) = 1$ and
$\|x - bzb\| \leq \epsilon$. It is sufficient to accomplish this
only for a dense set of elements $z$, so let $z = \sum_g x_gT_g$
where each $x_g$ is an element of $\A$ and $x_g \neq 0$ for only
finitely many $g$. (Such sums are clearly dense in $\A \times_r G$.)
\medskip

Let $(a_\lambda)$ be a decreasing net in $\A$ which excises $f$ and
satisfies $f(a_\lambda) = 1$ for all $\lambda$ [1, Proposition 2.2].
We claim that $\|f(x_e)a_\lambda^2 - a_\lambda za_\lambda\| \to 0$, which
will complete the proof. Observe first that
$\|f(x_e)a_\lambda^2 - a_\lambda x_ea_\lambda\| \to 0$ since $(a_\lambda)$
excises $f$. Since
$z$ is a finite sum and $x_e = x_eT_e$ we now need only show that
$\|a_\lambda (x_gT_g)a_\lambda\| \to 0$ for each $g \neq e$. But
$a_\lambda (x_gT_g)a_\lambda = (a_\lambda x_gb_\lambda)T_g$ where $b_\lambda =
\theta_g(a_\lambda)$ and $(b_\lambda)$ is a decreasing net which excises
$f\circ\theta_{g^{-1}}$. By hypothesis, $f$ and $f\circ\theta_{g^{-1}}$
are inequivalent, so Lemma 1 now implies $\|a_\lambda x_gb_\lambda\| \to 0$.
Thus $\|(a_\lambda x_gb_\lambda)T_g\| \to 0$, as desired.\hfill\hal
\bigskip

Note that in the proof of the reverse direction of Theorem 2, the fact
that $z = \sum x_gT_g$ can be ``compressed'' to $x_e$ implies that the
unique extension $f'$ of $f$ satisfies $f'(z) = f(x_e)$.
\medskip

By [12], the pure state space of any simple, separable C*-algebra $\A$ is
homogeneous --- given any two inequivalent pure states $f$ and $g$, there
is an automorphism $\omega$ of $\A$ such that $f = g\circ\omega$. In the
following corollary we need an even more powerful ``strong transitivity''
result which was proven in [7, Theorem 7.5] for a certain class of simple,
separable C*-algebras. By combining the methods of the two papers, that
result can be achieved for all simple, separable C*-algebras [11].
\medskip

The result we need states the following: {\it Suppose $\A$ is a simple,
separable C*-algebra and $(\pi_n)$ and $(\rho_n)$ are sequences of
irreducible representations such that the $\pi_n$ are mutually
inequivalent, as are the $\rho_n$. Then there is an automorphism
$\omega$ of $\A$ such that $\pi_n$ is equivalent to $\rho_n\circ\omega$
for all $n$.} It can be proven by first replacing every result in [12]
involving a single pure state $f$ (or a pair of pure states $f$ and $g$)
with a corresponding result involving a finite set of mutually inequivalent
pure states $f_i$ (or a pair of such sets, with, for each $i$, $f_i$ and
$g_i$ related in the same way that $f$ and $g$ were). The only difference
in the proofs will be that every application of Kadison's transitivity
theorem will now require Sakai's version for a finite family of mutually
inequivalent pure states [20, Theorem 1.21.16]. This achieves a proof
of [7, Theorem 7.3] for any simple separable C*-algebra, which can
be converted into a proof of the result we need in the same way that
this is done in [7, Theorem 7.5].
\bigskip

\noindent {\bf Corollary 3.} {\it Let $\A$ be a simple, separable, unital
C*-algebra and let $f$ and $g$ be inequivalent pure states on $\A$. Then
there is a simple, separable, unital C*-algebra $\B$ which unitally
contains $\A$ such that $f$ and $g$ have unique state extensions to $\B$
and these extensions are equivalent.}
\medskip

\noindent {\it Proof.} Observe first that if $f_1$ and $f_2$ are
equivalent pure states on $\A$ and $f_1$ has a unique state extension
to $\B$, then so does $f_2$. Indeed, if $u \in \A$ is a unitary
such that $f_1 = u^*f_2u$ then $u$ also
belongs to $\B$ and conjugates the set of extensions of $f_1$ with
the set of extensions of $f_2$.
\medskip

It follows from the result quoted before this corollary that given any
sequence $(\pi_n)$ ($n \in \Z$) of inequivalent irreps, there is an
automorphism $\omega$ of $\A$ such that $\pi_n \circ \omega$ is
equivalent to $\pi_{n+1}$ for all $n$. Now since $\A$ is simple and
separable and has inequivalent pure states it cannot be type I, and
therefore it has uncountably many inequivalent irreps [16, Corollary 6.8.5].
Let $(\pi_n)$ be any sequence of mutually inequivalent irreps such that
$\pi_1$ and $\pi_2$ are the GNS irreps arising from $f$ and $g$, and
find an automorphism $\omega$ as above. Then $g$ is equivalent to
$f\circ \omega$, and neither $f$ nor $g$ is equivalent to itself
composed with any nonzero power of $\omega$.
\medskip

Define $\theta: \Z \to {\rm Aut}(\A)$ by $\theta_n = \omega^n$ and let
$\B_0$ be the crossed product of $\A$ by this action of $\Z$. By
Theorem 2, $f$ and $g$ extend uniquely to $\B_0$, and therefore so do
all pure states equivalent to $g$, in particular $f\circ\omega$, by the
comment at the start of this proof. Now if $f'$ and $(f\circ\omega)'$
are the unique extensions of $f$ and $f\circ\omega$ to $\B_0$ then
$$\eqalign{f'\left(T_1\left(\sum x_nT_n\right)T_{-1}\right)
&= f'\left(\sum \omega(x_n)T_n \right)\cr
&= f(\omega(x_0)) = (f\circ\omega)'\left(\sum x_nT_n\right)\cr}$$
for any finite sum $\sum x_nT_n$, using the comment following Theorem 2.
This shows that $T_1f'T_{-1} = (f\circ\omega)'$, and hence $f'$ is
equivalent to $(f\circ\omega)'$.  Using the first paragraph again,
we see that the unique extensions of $f$ and $g$ to $\B_0$ are
equivalent.
\medskip

$\B_0$ is clearly separable and unital, and it unitally contains $\A$,
so we have proven the corollary except that $\B_0$ might not be simple.
To handle this, let $\B = \B_0/\I$ where $\I$ is any maximal ideal. Then
$\B$ is simple, and because $\A$ is simple (and contains the unit) the
quotient map takes it isomorphically into $\B$. Since every state on $\B$
becomes a state on $\B_0$ when composed with the quotient map, it follows
that any pure state on $\A$ which extends uniquely to $\B_0$ must extend
uniquely to $\B$. Also, since $f'$ and $(f \circ \omega)'$ are equivalent
via the unitary $T_1$, it follows that the extensions of $f$
and $f\circ\omega$ to $\B$ are equivalent via the unitary $T_1 + \I$,
which implies that the extensions of $f$ and $g$ to $\B$ are equivalent.
So $\B$ has all claimed properties.\hfill\hal
\bigskip
\bigskip

\noindent {\bf 2. The counterexample}
\bigskip

A subset $S$ of $\aleph_1$ is said to be {\it closed} if for
every countable $S_0 \subset S$ we have $\sup S_0 \in S$. It
is {\it unbounded} if for every $\alpha \in \aleph_1$ there
exists $\beta \in S$ such that $\beta > \alpha$.
\medskip

We call a nested transfinite sequence of C*-algebras $(\A_\alpha)$
{\it continuous} if for every limit ordinal $\alpha$ we have
$\A_\alpha = \overline{\bigcup_{\beta < \alpha} \A_\beta}$.
\bigskip

\noindent {\bf Lemma 4.} {\it Let $(\A_\alpha)$, $\alpha < \aleph_1$, be
a continuous nested transfinite sequence of separable C*-algebras and set
$\A = \bigcup \A_\alpha$. Then $\A$ is a C*-algebra, and if $f$ is
a pure state on $\A$ then $\{\alpha:
f$ restricts to a pure state on $\A_\alpha\}$ is closed and unbounded.}
\medskip

\noindent {\it Proof.} We observe first that $\A$ is automatically complete.
For any sequence $(x_n) \subset \A$ we can find indices $\alpha_n$
such that $x_n \in \A_{\alpha_n}$, and then $(x_n) \subset \A_\alpha$
for $\alpha = \sup \alpha_n$. Thus if $(x_n)$ is Cauchy, its limit
belongs to $\A_\alpha$ and hence to $\A$. This shows that $\A$ is a
C*-algebra.
\medskip

Let $f$ be a pure state on $\A$ and
let $S = \{\alpha: f$ restricts to a pure state on $\A_\alpha\}$. First
we verify that $S$ is closed. Suppose $S_0 \subset S$ is countable and
let $\alpha = \sup S_0$. If $f|_{\A_\alpha}$ is not pure then we
can write $f|_{\A_\alpha} = (f_1 + f_2)/2$ where $f_1$
and $f_2$ are distinct states on $\A_\alpha$. Now
$\bigcup_{\beta \in S_0} \A_\beta$ is dense in $\A_\alpha$ by
continuity (unless $\alpha \in S_0$, when this is true vacuously).
Thus there exists $\beta \in S_0$ such that $f_1|_{\A_\beta} \neq
f_2|_{\A_\beta}$. But then $f|_{\A_\beta} = (f_1|_{\A_\beta}
+ f_2|_{\A_\beta})/2$ contradicts the fact that $f|_{\A_\beta}$
is pure. We conclude that $f|_{\A_\alpha}$ is pure and hence that
$S$ is closed.
\medskip

Next observe that there is a sequence $(x_n)$ in $\A$ such
that $|f(x_n)|/\|x_n\| \to 1$. Then $(x_n) \subset \A_\alpha$ for some
$\alpha < \aleph_1$, so that the restriction of $f$ to $\A_\beta$
has norm one, and hence is a state, for all $\beta \geq \alpha$.
Thus, without loss of generality, in proving unboundedness we can assume
the restriction of $f$ to any $\A_\alpha$ is a state.
\medskip

Let $\alpha < \aleph_1$. We first claim that for any $x \in \A_\alpha$,
for sufficiently large $\beta \geq \alpha$ we have
$$f|_{\A_\beta} = (f_1 + f_2)/2 \qquad \Rightarrow \qquad f_1(x) = f_2(x)$$
whenever $f_1$ and $f_2$ are states on $\A_\beta$.
That is, for each $x \in \A_\alpha$ there exists $\alpha' \geq \alpha$
such that the above holds for all $\beta \geq \alpha'$. If the displayed
condition holds for all states $f_1$ and $f_2$ on $\A_\beta$ then we say
that $f|_{\A_\beta}$ is {\it pure on $x$}.
\medskip

Suppose the claim fails for some $x \in \A_\alpha$. Then there exist
$\epsilon > 0$ and an unbounded set $T \subset \aleph_1$ together with
states $f_1^\beta$ and $f_2^\beta$ on $\A_\beta$ for all $\beta \in T$,
such that
$$f|_{\A_\beta} = (f_1^\beta + f_2^\beta)/2\qquad{\rm and}
\qquad |f_1^\beta(x) - f(x)| \geq \epsilon.\eqno{(*)}$$
(If no such $\epsilon$ and $T$ existed, then for each $n \in \N$ we
could find $\alpha_n < \aleph_1$ such that for any $\beta \geq \alpha_n$,
no states $f_1^\beta$ and $f_2^\beta$ on $\A_\beta$
satisfy ($*$) with $\epsilon = 1/n$. Then $f|_{\A_\beta}$ would be
pure on $x$ for all $\beta \geq \sup \alpha_n$, contradicting our
assumption that the claim fails for $x$.) Now let $\U$ be an
ultrafilter on $T$ which contains the set $\{\beta \in T:
\beta > \beta_0\}$ for each $\beta_0 < \aleph_1$, and define states
$g_1$ and $g_2$ on $\A$ by $g_1 = \lim_\U f_1^\beta$
and $g_2 = \lim_\U f_2^\beta$, where the limits are taken
pointwise on elements of $\A$. Then $f = (g_1 + g_2)/2$
and $|g_1(x) - f(x)| \geq \epsilon$, hence $f \neq g_1$,
contradicting purity of $f$. (This argument could also be carried
out using universal nets.) This establishes the claim.
\medskip

Now for any $\alpha < \aleph_1$, since $\A_\alpha$ is separable we can
find a dense sequence $(x_n) \subset \A_\alpha$, and the claim implies
that for sufficiently large $\beta$, $f|_{\A_\beta}$ is pure on every
$x_n$. By density, $f|_{\A_\beta}$ is then pure on every $x \in \A_\alpha$.
Let $\alpha^*$ be the least ordinal larger than $\alpha$ such that
$\beta \geq \alpha^*$ implies that $f|_{\A_\beta}$ is pure on every
$x \in \A_\alpha$.
\medskip

Finally, to see that $S$ is unbounded, fix $\alpha < \aleph_1$.
Let $\alpha_1 = \alpha^*$, $\alpha_2 = \alpha^{**}$, etc., and
let $\beta = \sup \alpha_n$. Then $f|_{\A_\beta}$ is pure on every
$x \in \bigcup \A_{\alpha_n}$, and hence by continuity $f|_{\A_\beta}$
is a pure state on $\A_\beta$. This shows that $\beta \in S$,
as desired.\hfill\hal
\bigskip

A subset of $\aleph_1$ is {\it stationary} if it intersects every
closed unbounded subset of $\aleph_1$. We require the following
version of the diamond principle ($\diamondsuit$): there exists a
transfinite sequence of functions $h_\alpha: \alpha \to \aleph_1$
($\alpha < \aleph_1$) such that for any function $h: \aleph_1
\to \aleph_1$ the set $\{\alpha: h|_\alpha = h_\alpha\}$ is
stationary. (See [9, (22.20)] or [13, Exercise II.51].)
\medskip

Let $S(\A)$ denote the set of states on a C*-algebra $\A$.
\bigskip

\noindent {\bf Theorem 5.} {\it Assume $\diamondsuit$. Then there
is a counterexample to Naimark's problem which is generated by
$\aleph_1$ elements.}
\medskip

\noindent {\it Proof.} Let $(h_\alpha)$ be a transfinite sequence
of functions which verifies $\diamondsuit$ in the form given above.
For $\alpha < \aleph_1$ we recursively construct a continuous nested
transfinite sequence of simple separable unital C*-algebras $\A_\alpha$,
all with the same unit;
pure states $f_\alpha$ on $\A_\alpha$ with the property that
$\alpha < \beta$ implies $f_\beta$ is the unique state extension
of $f_\alpha$; and injective functions $\phi_\alpha: S(\A_\alpha)
\to \aleph_1$ as follows. Let $\A_0$ be any simple, separable,
infinite dimensional, unital C*-algebra and let $f_0$ be any pure
state on $\A_0$. Since $\diamondsuit$ implies the continuum hypothesis
(see [9] or [13]),
$S(\A_0)$ has cardinality at most $\aleph_1$, so there exists an
injective function from $S(\A_0)$ into $\aleph_1$; let $\phi_0$ be
any such function.
\medskip

To proceed from stage $\alpha$ of the construction to stage $\alpha + 1$
when $\alpha$ is a limit ordinal, first check whether there is a
pure state $g_\alpha$ on $\A_\alpha$, not equivalent to $f_\alpha$,
such that $h_\alpha(\beta) = \phi_\beta(g_\alpha|_{\A_\beta})$ for all
$\beta < \alpha$. (By injectivity of all $\phi_\beta$, there is at most
one such $g_\alpha$.) If so, let $\A_{\alpha + 1}$ be the C*-algebra
$\B$ given by Corollary 3 with $f = f_\alpha$ and $g = g_\alpha$.
Let $f_{\alpha + 1}$ be the unique extension of $f_\alpha$ to
$\A_{\alpha + 1}$ and as above let $\phi_{\alpha + 1}$ be any
injective function from $S(\A_{\alpha + 1})$ into $\aleph_1$.
If there is no such state $g_\alpha$, and whenever $\alpha$ is a
successor ordinal (or $\alpha = 0$), let $\A_{\alpha + 1}
= \A_\alpha$, $f_{\alpha + 1} = f_\alpha$, and $\phi_{\alpha + 1}
= \phi_\alpha$. At limit ordinals $\alpha$, let $\A_\alpha =
\overline{\bigcup_{\beta < \alpha} \A_\beta}$ (it is standard that
$\A_\alpha$ will be simple, given that each $\A_\beta$ is simple),
define $f_\alpha$ by requiring $f_\alpha|_{\A_\beta} = f_\beta$ for
all $\beta < \alpha$ ($f_\alpha$ will be pure, by the argument in the
second paragraph of the proof of Lemma 4), and as before let
$\phi_\alpha$ be any injective function from $S(\A_\alpha)$ into
$\aleph_1$.
\medskip

Let $\A = \bigcup_{\alpha < \aleph_1} \A_\alpha$ and define $f \in S(\A)$
by $f|_{\A_\alpha} = f_\alpha$. Here $f$ is well-defined since $f_\beta$
is an extension of $f_\alpha$ whenever $\alpha < \beta$. Also, $f$ is a
pure state, again by the reasoning in the second paragraph of the proof
of Lemma 4 plus the fact that each $f_\alpha$ is pure. We claim
that every pure state $g$ on $\A$ is equivalent to $f$. To see this,
it is enough to verify that $g|_{\A_\alpha}$ is equivalent to
$f_\alpha$ for some $\alpha$; then since $f_\alpha$ extends uniquely
to $\A$ the same must be true of $g|_{\A_\alpha}$ (see the first
paragraph of the proof of Corollary 3), and the unitary in $\A_\alpha$
which implements the equivalence of $f_\alpha$ and $g|_{\A_\alpha}$
must then also implement an equivalence between $f$ and $g$. Now
define $h: \aleph_1 \to \aleph_1$ by setting $h(\alpha) =
\phi_\alpha(g|_{\A_\alpha})$. Let $S$ be the set of limit ordinals
$\alpha$ such that $g|_{\A_\alpha}$ is pure. According to Lemma 4,
$S$ is closed and unbounded. (The intersection of any two closed
unbounded sets is always closed and unbounded.) Therefore, by
$\diamondsuit$ there exists a limit ordinal $\alpha$ such that
$g|_{\A_\alpha}$ is pure and $h_\alpha = h|_\alpha$, i.e.,
$$h_\alpha(\beta) = \phi_\beta(g|_{\A_\beta})$$
for all $\beta < \alpha$. If $g|_{\A_\alpha}$ is equivalent to
$f_\alpha$ then we are done, and otherwise the construction at
stage $\alpha$ guarantees that $f_{\alpha + 1}$ is equivalent to
the unique extension of $g|_{\A_\alpha}$ to $\A_{\alpha + 1}$,
which must be $g|_{\A_{\alpha + 1}}$. This completes the proof
that $f$ and $g$ are equivalent.
\medskip

Finally, $\A$ is infinite dimensional and unital, so it cannot
be isomorphic to any $\K(H)$.\hfill\hal
\bigskip

We conclude by observing that if the continuum hypothesis fails then
there is no counterexample to Naimark's problem which is generated by
$\aleph_1$ elements.
\bigskip

\noindent {\bf Proposition 6.} {\it Let $\A$ be a counterexample
to Naimark's problem. Then $\A$ cannot be generated by fewer than
$2^{\aleph_0}$ elements.}
\medskip

\noindent {\it Proof.} Since $\A$ cannot be type I, it follows
that there is a subalgebra $\B$ of $\A$ and a surjective
$*$-homomorphism of $\B$ onto the Fermion algebra $\F =
\bigotimes_1^\infty M_2(\C)$ [16, Corollary 6.7.4]. Let $f_1$
and $f_2$ be the pure states on $M_2(\C)$ defined by
$$f_1\left(\left[\matrix{a_{11}&a_{12}\cr a_{21}&a_{22}\cr}\right]\right)
= a_{11}
\qquad{\rm and}\qquad
f_2\left(\left[\matrix{a_{11}&a_{12}\cr a_{21}&a_{22}\cr}\right]\right)
= a_{22}.$$
Then for any sequence of indices $i_n \in \{1, 2\}$ the product state
$f = \bigotimes_n f_{i_n}$ is a pure state on $\F$. Moreover,
$\|f - f'\| = 2$ for any two distinct states of this form. Lifting to
$\B$ and extending to $\A$, we obtain a family of
$2^{\aleph_0}$ pure states on $\A$, the distance between any two of
which is 2. As $\A$ is a counterexample to Naimark's problem, there
exists a family of $2^{\aleph_0}$ unitaries in the unitization $\tilde{\A}$
of $\A$ which conjugate this family of pure states to a single (arbitrary)
pure state $g$ on $\A$. Uniform separation of the pure states implies
uniform separation of the conjugating unitaries: if $g = u_1^*f_1u_1
= u_2^*f_2u_2$ then for all $x \in \A$
$$|f_1(x) - f_2(x)| = |f_1(x - u^*xu)| \leq
|f_1(x - u^*x)| + |f_1(u^*x - u^*xu)| \leq 2\|1_\A - u\|\|x\|$$
where $u = u_2^*u_1$, so $\|f_1 - f_2\| \leq 2\|1_\A - u\| =
2\|u_1 - u_2\|$. Thus $\tilde{\A}$, and hence $\A$, cannot contain a
dense set with fewer than $2^{\aleph_0}$ elements, and consequently
it cannot be generated by fewer than $2^{\aleph_0}$ elements.\hfill\hal
\bigskip

\noindent {\bf Corollary 7.} {\it If the continuum hypothesis fails,
then no C*-algebra generated by $\aleph_1$ elements is a counterexample
to Naimark's problem. The existence of a counterexample to Naimark's
problem which is generated by $\aleph_1$ elements is independent of ZFC.}
\bigskip

The first assertion of Corollary 7 follows immediately from
Proposition 6, and the second assertion follows from Theorem 5 plus
the first assertion, together with the fact that both diamond and
the negation of the continuum hypothesis are consistent with ZFC
(assuming ZFC is consistent).
\bigskip
\bigskip
%##

\noindent [1] C.\ Akemann, J.\ Anderson, and G.\ K.\ Pedersen,
Excising states of C*-algebras, {\it Canad.\ J.\ Math.\ \bf 38} (1986),
1239-1260.
\medskip

\noindent [2] J.\ Anderson, Extensions, restrictions, and
representations of states on C*-algebras, {\it Trans.\ Amer.\
Math.\ Soc.\ \bf 249} (1979), 303-329.
\medskip

\noindent [3] W.\ Arveson, {\it An Invitation to C*-Algebras},
Springer-Verlag GTM {\bf 39}, 1976.
\medskip

\noindent [4] K.\ J.\ Devlin, {\it The Axiom of Constructibility:
A Guide for the Mathematician}, Springer-Verlag LNM {\bf 617}, 1977.
\medskip

\noindent [5] J.\ Dixmier, Sur les C*-alg\`ebres, {\it Bull.\
Soc.\ Math.\ France \bf 88} (1960), 95-112.
\medskip

\noindent [6] J.\ M.\ G.\ Fell, C*-algebras with smooth dual,
{\it Illinois J.\ Math.\ \bf 4} (1960), 221-230.
\medskip

\noindent [7] H.\ Futamura, N.\ Kataoka, and A.\ Kishimoto,
Homogeneity of the pure state space for separable C*-algebras,
{\it Internat.\ J.\ Math.\ \bf 12} (2001), 813-845.
\medskip

\noindent [8] J.\ Glimm, Type I C*-algebras, {\it Ann.\ of Math.\
\bf 73} (1961), 572-612.
\medskip

\noindent [9] T.\ Jech, {\it Set Theory}, Academic Press, 1978.
\medskip

\noindent [10] I.\ Kaplansky, The structure of certain operator
algebras, {\it Trans.\ Amer.\ Math.\ Soc.\ \bf 70} (1951), 219-255.
\medskip

\noindent [11] A.\ Kishimoto, personal communication.
\medskip

\noindent [12] A.\ Kishimoto, N.\ Ozawa, and S.\ Sakai,
Homogeneity of the pure state space of a separable C*-algebra,
{\it Canad.\ Math.\ Bull.\ \bf 46} (2003), 365-372.
\medskip

\noindent [13] K.\ Kunen, {\it Set Theory: An Introduction to
Independence Proofs}, North-Holland, 1980.
\medskip

\noindent [14] M.\ A.\ Naimark, Rings with involutions, {\it Uspehi
Matem.\ Nauk \bf 3} (1948), 52-145 (Russian).
\medskip

\noindent [15] {---------}, On a problem in the theory of rings with
involution, {\it Uspehi Matem.\ Nauk \bf 6} (1951), 160-164 (Russian).
\medskip

\noindent [16] G.\ K.\ Pedersen, {\it C*-Algebras and Their Automorphism
Groups}, Academic Press, 1979.
\medskip

\noindent [17] A.\ Rosenberg, The number of irreducible representations
of simple rings with no minimal ideals, {\it Amer.\ J.\ Math.\ \bf 75}
(1953), 523-530.
\medskip

\noindent [18] S.\ Sakai, A characterization of type I C*-algebras,
{\it Bull.\ Amer.\ Math.\ Soc.\ \bf 72} (1966), 508-512.
\medskip

\noindent [19] {---------}, On type I C*-algebras, {\it Proc.\
Amer.\ Math.\ Soc.\ \bf 18} (1967), 861-863.
\medskip

\noindent [20] {---------}, {\it C*-Algebras and W*-Algebras},
Springer-Verlag, 1971.
\medskip

\noindent [21] S.\ Shelah, Infinite abelian groups, Whitehead
problem and some constructions, {\it Israel J.\ Math.\ \bf 18} (1974),
243-256.
\medskip

\noindent [22] {---------}, Uncountable constructions for B.A.,
e.c.\ groups and Banach spaces, {\it Israel J.\ Math.\ \bf 51} (1985),
273-297.
\bigskip
\bigskip

\noindent Department of Mathematics

\noindent UCSB

\noindent Santa Barbara, CA 93106 USA

\noindent akemann@math.ucsb.edu
\bigskip

\noindent Department of Mathematics

\noindent Washington University

\noindent St.\ Louis, MO 63130 USA

\noindent nweaver@math.wustl.edu
\end